\documentclass[10pt,a4paper]{article}
\usepackage{amsmath,amssymb}
\begin{document}
\setlength{\baselineskip}{1.5em}
\begin{center} 
{\bfseries The Complex of Words and Nakaoka Stability}\\
Moritz C. Kerz\\
\textit{mokerz@students.uni-mainz.de}
\end{center}
\newtheorem{theo}{Theorem}
\newtheorem{lemma}{Lemma}
\newtheorem{defi}{Definition}
\newtheorem{coro}{Corollary}
\renewcommand{\labelenumi}{(\roman{enumi})}
{\bfseries Abstract}\\
We give a new simple proof of the exactness of the complex of injective words and use it to prove
Nakaoka's homology stability for symmetric groups. The methods are generalized to show acyclicity in low 
degrees for the complex of words in ''general position''.\\

\textbf{1. An elementary proof of Nakaoka stability}\\

We will be concerned with the complex $C_*(m)$, where $C_n(m)$ is the free abelian group generated by elements
 $(x_1,\ldots,x_n)$, where the $x_i$ are pairwise different natural numbers from $1$ to $m$. $C_0(m) = \mathbb{Z}$.
 The differential is given by
\[
d (x_1,\ldots, x_n) = \sum_{j=1}^{n}(-1)^{j+1}(x_1,\ldots, x_{j-1},x_{j+1},\ldots, x_n)\: .
\]
In discrete mathematics the following theorem is well established as the homology of $C_*(M)$ is given by the simplicial homology of the 
(shellable) poset of injective words. Shellability reduces the simplicial homology groups to those of a wedge of m-spheres. 
 We refer to [1], [3], [7] for exact definitions and statements.\\
Our proof of Theorem 1 is new and rather straightforward compared to Farmer's original elementary proof. 

\begin{theo}{(F.D.~Farmer [3])}
The homology of $C_{*}(m)$ vanishes except in degree $m$.
\end{theo} 

We have to introduce some notations which will be used throughout the paper.\\
A chain c is called a term, if there exists an $N\in \mathbb{Z}$ and $x_1,\ldots, x_n \in \{1,\ldots, m\}$ with $c = N\, (x_1,\ldots,x_n)$.\\
As $C_*(m)$ has a canonical basis all our sum decompostions correspond to partitions of the basis. We also speak about the appearance of numbers in a chain. For example
$2$ appears in $(2,3)+4\, (5,1)$ but $4$ does not. \\
Although our complex $C_*(m)$ has no obvious cup product, we have a partially defined product. If $c\in C_n(m)$, $c'\in C_l(m)$ are elementary $c=N \,(x_1,\ldots x_n)$,
 $c= M\, (x'_1,\ldots ,x'_l)$, we define $c\, c' :=N\, M\, (x_1,\ldots, x_n, x'_1,\ldots x'_l)$ if the numbers $x_1,\ldots, x_n, x'_1,\ldots x'_l$ are paiwise different. 
This constuction extends bilinearly to arbitrary chains $c$, $c'$ for which the numbers appearing in both of them are distinct. \\
There is a Leibniz rule for $c\in C_n(m)$ and $c'\in C_l(m)$
\[
d(c\, c') = d(c)\, c' + (-1)^n c\, d(c')\: .
\]

\vspace{0.3cm}

Proof of Theorem 1.
The exactness in degree $0$ is clear. For the rest we use induction on $m$. For the case $m=2$ we have to check that 
\[
C_2(2)\longrightarrow C_1(2)\longrightarrow \mathbb{Z}
\] 
is exact. But $\mbox{ker}(d:C_1(2)\rightarrow \mathbb{Z})$ is generated by elements of the form
\[
(x)-(x') = d\, (x,x')
\]
with $x\ne x'$.

For the induction step we use a straightforward lemma:

\begin{lemma}
If we have a number $x\in\{1,\ldots,m\}$ which does not appear in a cycle $c\in C_n(m)$, it is a boundary.
\end{lemma}
Proof. According to the Leibniz rule $c =c + (x)\, d(c)= d( (x)\, c)$, since 
$d(c) = 0$.\hfill $\Box$\\

Given an arbitrary cycle $c$ of degree $n<m$ we have to show that in order to apply the lemma we can eliminate a number from $c$ by adding boundaries. Therefore we will push a fixed number $x\in\{1,\ldots,m\}$ to the right until it vanishes from the cycle.\\
If $x$ appears somewhere in the cycle at the first index, write
\[
c=\sum_{j}(x)\, c_j + c'
\]
Where the $c_j$ are elementary and $c'$ does not have $x$ at the first index. To each $c_j$ choose a number $x_j\in\{1,\ldots,m\}$, $x_j\ne x$ which does not appear in $c_j$.
Then
\[
c - d(\sum_{j}(x_j)\, (x)\, c_j)  =  c' + \sum_j (x_j)\ c_j - (x_j)\, (x)\, dc_j \; .
\] 
Now $x$ does not appear at the first index anymore. The next steps until the vanishing of $x$ are similar. Suppose $x$ does not appear at the first $i>0$ indices of $c$. Now we can write
\[
c = \sum_j s_j\, (x)\, c_j + c'\: .
\]

Where the $c_j$ are elementary and different, the $s_j$ have length $i$, $x$ does not appear in $s_j$ and $x$ does not appear at the first $i+1$ entries of $c'$.
One calculates:
\[
0=d\, c = \sum_j\left[ (d\, s_j)\, (x)\, c_j +(-1)^{i} s_j\, c_j +(-1)^{i+1} s_j\, (x)\, d\, c_j \right] +d\, c'\; .
\]
If we forget those terms for which $x$ does not appear at the $i$-th index, this equation implies $d\, s_j = 0$ for all $j$. Since $length(s_j)=i<m-length((x)\, c_j)$,
 there are by induction $s_j'$ with $d\, s_j' = s_j$ and such that the following products make sense.
\[
z:=\sum_j s_j'\, (x)\, c_j
\]
In the cycle
\[
c - d\, z = c' - \sum_j \left[(-1)^{i+1}s_j'\, c_j + (-1)^i s_j'\, (x)\, d\, c_j\right]
\]
$x$ does not appear at the first $i+1$ indices.
\hfill $\Box$\\

Using the corresponding two hyperhomology spectral sequences for the natural action of the symmetric group $\Sigma_m$ on our complex $C_*(m)$ (cf [2]) one can now obtain a stability result due to Nakaoka.
A similar proof was given by Maazen (thesis).

\begin{theo}{(Nakaoka [5])}
$H_m(\Sigma_{n-1})=H_m(\Sigma_n)$ for $m<n/2$.
\end{theo}
Proof. We use induction on $n$ . It is well known for $n=3$ 
\[
H_1(\Sigma_2) = H_1(\Sigma_3) = \mathbb{Z}/(2) \: .
\]
For $n\ge 4$ define $C'_l(n):=C_{l+1}(n)$ for $l\ge 0$. Then
\[
H_m(\Sigma_n,\mathbb{Z})=H_m(\Sigma_n,C'_*(n))
\] 
when $m<n-1$ because of Theorem 1. The other spectral sequence of the bi-complex gives for $E_{*,*}^1$:

\[
\begin{array}{ccccc}
 &  \cdots & & & \\
H_2(\Sigma_n,C'_0) & H_2(\Sigma_n,C'_1) & \cdots & H_2(\Sigma_n,C'_{n-2}) & H_2(\Sigma_n,C'_{n-1})\\
H_1(\Sigma_n,C'_0) & H_1(\Sigma_n,C'_1) & \cdots & H_1(\Sigma_n,C'_{n-2}) & H_1(\Sigma_n,C'_{n-1})\\
H_0(\Sigma_n,C'_0) & H_0(\Sigma_n,C'_1) & \cdots & H_0(\Sigma_n,C'_{n-2}) & H_0(\Sigma_n,C'_{n-1})

\end{array}
\]

Using Shapiro's Lemma we get:
\[
\begin{array}{ccccc}
 &  \cdots & & & \\
H_2(\Sigma_{n-1}) & H_2(\Sigma_{n-2}) & \cdots & H_2(\Sigma_{1}) & 0\\
H_1(\Sigma_{n-1}) & H_1(\Sigma_{n-2}) & \cdots & H_1(\Sigma_{1}) & 0\\
H_0(\Sigma_{n-1}) & H_0(\Sigma_{n-2}) & \cdots & H_0(\Sigma_1) & \mathbb{Z}
\end{array}
\]
The horizontal arrows can be computed as $0,1,0,1,\cdots$, since they are the sums of the signs in $d'_1,d'_2,\cdots$.\\
We have $E^2_{i,0} = H_i(\Sigma_{n-1})$ for $i\ge 0$ and by induction $E^2_{i,j}=0$ for $i<\frac{n-j-1}{2},\: 0<j<n-1$.\\
Finally, parts of the spectral sequence degenerate; $E^{\infty}_{i,0}=H_i(\Sigma_{n-1})$ for $n-i>2$,  such that $H_i(\Sigma_n)=H_i(\Sigma_{n-1})$ for $n/2>i$.
\hfill $\Box$

\vspace{0.4cm}

\textbf{2. Words in general position}

\vspace{0.3cm}
\renewcommand{\labelenumi}{(\roman{enumi})}

Let $\mathbf{n}:=\{i|1\le i \le n \}$.
We associate to every (nonempty) set $X$ a complex $F(X)_*$, the so-called
complex of words with alphabet $X$, as follows:
\begin{eqnarray*}
(F(X))_n & =  & \mathbb{Z}< \{f:\mathbf{n}\to X \}>\\
 d_n (f) & = & \sum_{k=1}^n (-1)^{n+1} f\circ \delta_k  \: .
\end{eqnarray*}
Here  $\delta_k:[n]\rightarrow [n+1]$ are the coface maps
\[
\delta_k(j)=\left\{ 
\begin{array}{ll} 
j\;\;& \mbox{if } 1\le j< k\\
j+1\;\; & \mbox{if } k\le j\
\end{array} 
\right.\: .
\]

It is immediate that the homology of $F(X)$ vanishes. For if $c\in F(X)_n$ is a cycle, $c = d((x_0)\, c)$ 
for arbitrary $x_0\in X$ (see Lemma 1).\\
For given $X$ certain subcomplexes of $F(X)$ can be used in hyperhomology spectral sequences as above, if their homology 
vanishes to some extent. These subgroups are determined by conditions which one could call ''general position conditions''.
For applications cf [8], [9].\\ 
It could be asked, how to explain the fact that the vanishing of homology is not affected by these conditions.
In order to give a general result we translate the proof of Theorem 1 into our more complicated setting.

\vspace{0.5cm}

{\bfseries Examples}\\
\begin{enumerate}
\item{Let $X$ be a finite set. The complex of injective words is $(F_{inj}(X))_n:=\{f\in F(X)_n | f\, \text{injective}\}$. 
According to Theorem 1 the homology of this subcomplex is zero except in degree $m=card(X)$ where 
\[\mbox{rank}\, (H_m(F(X)_*))= (-1)^{m}(1-\sum_{i=0}^{m-1}(-1)^i m (m-1)\cdots (m-i) )\: .  \]
}
This equals the number of fixed point free permutations of $X$.
\item{Let $F$ be a field and $V$ an $F$-vector space. 
The complex of vectors in general position is $(F_{gen}(V))_n:=\{f\in F(V)_n |f \mbox{ in general position}\}$ .\\
If $F$ is infinite the complex is exact. A general vanishing result is contained in our main theorem.}
\end{enumerate}
We introduce axioms for elements of a given (nonempty) set $X$ to be in general position.

\begin{defi}[General position] 
Let $G_{l,m}$ be relations in $l+m$ variables from $X$; $l> 0$, $m\ge 0$.
The family $(G_{l,m})_{l,m\in \mathbb{N}}$ is called a general position relation if the following
properties are satisfied.

Let $x,y,z$ be finite sequences of elements of $X$ of length $l,m,n$:
\begin{enumerate}
\item{$G_{n,m}$ is symmetric in the first $n$ and last $m$ arguments.}
\item{If $G_{l+m,n}(x,y;z)$ then $G_{l,m+n}(x;y,z)$. If $G_{l,m+n}(x;y,z)$ then $G_{l,n}(x;z)$.}
\item{If $G_{l,m+n}(x;y,z)$ and $G_{m,n}(y;z)$ then $G_{l+m,n}(x,y;z)$.}
\end{enumerate}
\end{defi}

For fixed $G$ we say $x$ is in $G$-general position to $y$ iff $G_{l,m}(x;y)$.\\
Given $x\in F(X)_l$ and $y\in F(X)_m$ we say $x$ is in $G$-general position to $y$ if every term of $x$ is in
$G$-general postition to every term in $y$ (for $l=0$ we demand nothing).
\begin{defi}
Given a finite sequence $b$ of elements of $X$ the corresponding complex of chains in general position to $b$ is 
$(F_{G}(X;b))_n:=\{f\in F(X)_n | f$ is in $G$-general position to $b \}\, $.
\end{defi}

Before we can state the main theorem we have to introduce an invariant which determines an upper bound for the 
vanishing of the homology of $F(X)_*$.

\begin{defi}
Given a general position relation $G$ we define $|G|$ to be the smallest natural number $n\ge 0$ such that there is
a sequence $x$ of elements of $X$ with $length(x)=n$ such that there is no further element in $X$ which is in general
position to x.
\end{defi}

{\bfseries Examples}\\
\begin{enumerate}
\item{Example (i) is induced by saying $x$ is in $G^{inj}$-general position to $y$ if the underlying sets of $x$ and $y$ are
disjoint and the entries of $x$ are pairwise different. We have $|G^{inj}|=card(X)$ and $F_{G^{inj}}(X)_* = 
(F_{inj}(X))_*$ }

\item{Example (ii) is induced by saying $(x_1,\ldots ,x_i)$ is in $G^{vec}$-general position to 
$(y_1,\ldots , y_j)$ if we have for all $a_k\in F$,
$k=1,\ldots,i+j$ and only $dim(V)$ many of them nonvanishing 
\[
a_1\, x_1 + \cdots + a_i\, x_i + a_{i+1}\, y_1 + \cdots + a_{i+j}\, y_j=0
\] 
implies $a_k=0$ for all $k\in \{1,\ldots,i\}$.\\
If $card(F)=\infty$ or $dim(V)=\infty$ then $|G^{vec}|=\infty$.\\
Unfortunately the exact value of $|G^{vec}|$ is not known for all finite dimensional vector spaces over
finite fields. }
\end{enumerate}

\begin{lemma}
(a) For $dim(V)\ge card(F)$ we have $|G^{vec}|=dim(V)+1$.\\
(b) For $dim(V)=2$ we have $|G^{vec}|=card(F)+1$.
\end{lemma}
Proof of (a). Let $n=dim(V)+1$  and $(e_i)_{1\le i \le n-1}$ be a basis of $V$.
We show $|G^{vec}|\ge n$ first. Otherwise we had a sequence $(x_i)_{1\le i\le n-1}$, $x_i\in V$, such that
there does not exist a vector $x\in V$ in $G^{vec}$-general position to $(x_i)$. This can only be true if $(x_i)$
is a basis of $V$. But if it is a basis the element $x_1+x_2+\cdots +x_{n-1}$ would be in general position to it. 
Contradiction.\\
Now the sequence $f=(e_1,\ldots ,e_{n-1},e_1+e_2+\cdots +e_{n-1})$ is maximal in the sense of Definition 3. To see this
we have to show there is no vector $x=a_1\, e_1 + \cdots + a_{n-1}\, e_{n-1}$, $a_i\in F$, in $G^{vec}$-general position to
$f$. If there is an index $j$ such that $a_j=0$, $x$ is obviously not in general postion to $f$. But if no coefficient in $x$
vanishes two of them have to be equal since $dim(V)\ge card(F)$. In case $a_1=a_2$ we can write 
\[
x=a_1\, (e_1+\cdots +e_{n-1}) + (a_3-a_1)\, e_3 + \cdots (a_{n-1}-a_1)\, e_{n-1}\: . 
\] 
This shows $f$ is in fact maximal.\\
Proof of (b). For $dim(V)=2$ we have $x\in V$ in $G^{vec}$-general position to $y\in V$ iff $x\ne 0\ne y$ and
$[x]\ne [y]\in \mathbb{P}(V)$. So $|Gen^{vec}|$ is simply the number of rational points in $\mathbb{P}(V)/F$.
  \hfill $\Box$

\begin{theo}
For a set $X$ with general position condition $G$ and a finite sequence $a=(a_1,\ldots, a_l)$ of elements of $X$ the 
corresponding homology groups $H_m(F_{G}(X;a))$ vanish for $m\le (|G|-l-1)/2$.
\end{theo}
In contrast to Theorem 1 one should notice that it is in general 
not possible to erase the factor $1/2$ from the inequality of the theorem.

Proof. Denote the degree by $m$. The exactness at $m=0$ is trivial. We proceed by induction on $m\ge 1$. Let $c$ be a cycle
in $F_m(X;a_1,\ldots, a_l)$. \\
Case $m=1$:\\
We can suppose $c=(x)-(x')$. 
Because $1\le (|G|-l-1)/2$ we have $l+2<|G|$ so that there exists $y\in X$ in $G$-general position to
$(x,x',a_1,\ldots, a_l)$.\\
According to Definition 1(iii) $(y,x,x')$ is in $G$-general position to $(a_1,\ldots,a_l)$ and it is allowed to
write $d((y)\, c)=c$.

Induction step:\\ 
Choose $x\in X$ in $G$-general position to $(a_1,\ldots,a_l)$. \\
The simplest case is $c$ in $G$-general position to
$(x,a_1,\ldots,a_l)$. Here we can apply a construction similar to Lemma 1; we have $c=d((x)\, c)$, 
since $(x)\, c$ is in $G$-general position to $(a_1,\ldots,a_l)$. 

We reduce to this case by changing $c$ by boundaries. To be more precise we introduce a number $I(c)\in \{0,\ldots ,m\}$ 
which is $m$ iff the above applies, that is $c$ in $G$-general position to $(x,a_1,\ldots,a_l)$.\\
By adding boundaries be will see that we can increase $I(c)$.
 
For $d\in F_G(X;a)_n$ we define $I(d)\in\{0,\ldots,n\}$ as the 
the greatest natural number $i\le n$ such that $\pi_i(v)$ is in $G$-general position to $(x, \pi'_i(v), a_1,\ldots,a_l)$ 
for any term $v$ of $d$ ($\pi_i$ denotes the projection to the first $i$ entries and 
$\pi'_i$ the projection to the last $n-i$ entries).

Reduction to $I(c)>0$:\\ 
Suppose $I(c)=0$.
Let $c=\sum_j x_j$ with $x_j$ elementary. 
Choose for every $j\:$ $y_j\in X$ with $y_j$ in $G$-general position to $(x,x_j,a_1,\ldots,a_l)$. This is possible since
$\text{length}(x,x_j,a_1,\ldots,a_l)=1+m+l<|G|$. 
Clearly $I(c-d( \sum_j  (y_j)\, x_j))>0$.

Now suppose $m>I(c)>0$: \\
Write 
\[
c=\sum_j s_j\, x_j + x'
\]
such that exactly those terms $v$ of $c$ for which $I(v)>I(c)$ are in $x'$, $\text{length}(s_j)=I(c)$ and all $x_j$ 
are elementary and pairwise different. 

\begin{lemma}
$d(s_j)=0$ for all $j\:$.
\end{lemma}
Proof. The terms $v$ of $d(c)$ such that $I(v)<I(c)$ are exactly the terms of $d(s_j)\, x_j$.
In order to see this, notice that for a term $v$ of $d(c)$ we have $I(v)=I(c)-1$ iff the I(c)-th entry of $v$
is not in $G$-general position to $(x,\pi'_{I(c)}(v) , a_1,\ldots ,a_l)$.\\
Now projecting the identity $0=d(c)$ to the terms $v$ with $I(v)<I(c)$ we get $0=\sum_j s_j\, x_j$. Lemma 3 is proven. 
\hfill $\Box$ 

According to our assumtion $s_j$ is in $G$-general position to $(x,x_j,a_1,\ldots,a_l)$. 
By induction on $m$ we know there exist $s_j'$ with $d(s_j')=s_j$ and $s_j'$ in $G$-general position to 
$(x,x_j,a_1,\ldots,a_l)$, 
because $2 i + [(m-i)+l+1]\le 2 m + l \le |G|-1$.

Now 
\begin{eqnarray*}
I(c-d(\sum_j s_j'\, x_j)) & = & I(x'+(-1)^{I(c)} s_j'\, d(x_j) ) \\
& >&I(c)\: . 
\end{eqnarray*}
This finishes the reduction to $I(c)=m$ and therefore the induction step for $\text{length}(c)=m$ is accomplished.
\hfill $\Box$ \\

Theorem 3 could be useful in the generalization of Suslin's $GL$-stability [8] to finite fields. A thorough 
treatment seems to indicate the following result:\\
(i) Given $m\ge 0$ and $n\ge m$ we have $H_m(GL_n(F))=H_m(GL_{n+1}(F))$ for almost all finite fields $F$.\\
(ii) For given $m\ge 0$ $H_m(GL_{m-1}(F))\rightarrow H_m(GL_m(F))$ is surjective for almost all finite fields $F$.\\
In fact (i) has been proven -- using other methods -- by Quillen [6]. He did even show, that the statement is true for all fields
with more then $2$ elements. Similar results with weaker bounds are due to Maazen and van der Kallen [4].

\vspace{0.3cm}

\textbf{Acknowledgement}\\
The above results were obtained in the Seminar M\"uller-Stach 2004 at Mainz. I am indebted to Stefan 
M\"uller-Stach, Volkmar Welker, W.~van der Kallen, Philippe Elbaz-Vincent and Oliver Petras for their support.\\

\vspace{0.4cm}

{\bfseries References\\}

[1] Bj\"orner, Anders; Wachs, Michelle 
\textsl{On lexicographically shellable posets. } 
Trans. Amer. Math. Soc. 277 (1983), no. 1, 323--341.

[2] Brown, Kenneth S.
 \textsl{Cohomology of groups.}
Corrected reprint of the 1982 original. Graduate Texts in Mathematics, 87. 
Springer-Verlag, New York, 1994 

[3] Farmer, Frank D. 
\textsl{Cellular homology for posets.} 
Math. Japon. 23 (1978/79), no. 6, 607--613.

[4] van der Kallen, Wilberd 
\textsl{Homology stability for linear groups.} 
Invent. Math. 60 (1980), no. 3, 269--295.

[5] Nakaoka, Minoru 
\textsl{Decomposition theorem for homology groups of symmetric groups. }
Ann. of Math. (2) 71 1960 16--42.

[6] Quillen, Daniel; 
unpublished.

[7] Stanley, Richard P.
\textsl{Combinatorics and commutative algebra. } 
Second edition. Progress in Mathematics, 41. 
Birkh\"auser Boston, Inc., Boston, MA, 1996.

[8] Suslin, A. A. 
\textsl{Homology of ${\rm GL}\sb{n}$, characteristic classes and Milnor $K$-theory. } 
Algebraic $K$-theory, number theory, geometry and analysis (Bielefeld, 1982), 357--375, 
Lecture Notes in Math., 1046, 
Springer, Berlin, 1984. 

[9] Suslin, A. A. 
\textsl{$K\sb 3$ of a field, and the Bloch group. } 
Proc. Steklov Inst. Math. 1991, no. 4, 217--239.

\end{document}